\documentclass[aos,nameyear,dvips]{arximspdf}
\usepackage{mathbh,mathrsfs}
\usepackage{graphicx}

\doi{10.1214/09-AOS775}
\volume{38}
\issue{6}
\pubyear{2010}
\firstpage{3445}
\lastpage{3457}

\makeatletter

\newtheorem{theorem}{Theorem}
\newtheorem{lemma}{Lemma}
\newproclaim{remark}{Remark}

\def\eqref#1{(\ref{#1})}
\makeatother

\begin{document}
\begin{frontmatter}

\title{On optimality of the Shiryaev--Roberts procedure for detecting
a change in distribution\thanksref{T1}}
\runtitle{Optimality of the Shiryaev--Roberts procedure}
\thankstext{T1}{Supported in part by the US Army Research Office MURI Grant
W911NF-06-1-0094 and NSF Grant
CCF-0830419 at
the University of Southern California, Los Angeles.}

\begin{aug}
\author{\fnms{Aleksey S.} \snm{Polunchenko}\ead[label=e1]{polunche@usc.edu}}
\and
\author{\fnms{Alexander G.} \snm{Tartakovsky}\corref{}\ead[label=e2]{tartakov@usc.edu}}

\affiliation{University of Southern California, Los Angeles}

\address{Department of Mathematics\\
University of Southern California\\
3620 S. Vermont Ave, KAP-416F\\
Los Angeles, California 90089-2532\\ USA\\
\printead{e1}\\
\phantom{E-mail: }\printead*{e2}}

\runauthor{A. S. Polunchenko and A. G. Tartakovsky}
\end{aug}

\received{\smonth{4} \syear{2009}}
\revised{\smonth{8} \syear{2009}}

%
\begin{abstract}
In 1985, for detecting a change in distribution, Pollak introduced a
specific minimax performance metric and a
randomized version of the Shiryaev--Roberts procedure where the zero
initial condition is replaced by a random
variable sampled from the quasi-stationary distribution of the
Shiryaev--Roberts statistic. Pollak proved that this
procedure is third-order asymptotically optimal as the mean time to
false alarm becomes large. The question of whether
Pollak's procedure is strictly minimax for any false alarm rate has
been open for more than two decades, and there were
several attempts to prove this strict optimality. In this paper, we
provide a counterexample which shows that Pollak's
procedure is not optimal and that there is a strictly optimal procedure
which is nothing but the Shiryaev--Roberts procedure
that starts with a specially designed deterministic point.\looseness=-1
\end{abstract}



\end{frontmatter}

\section{Introduction and preliminaries}\label{s:Intro}

Changepoint problems deal with detecting changes in distributions of
observed data that occur at unknown points in time.
Let $X_1, X_2, \dots$ be the series of observations being monitored,
and let $\nu$ be the
serial number of the last pre-change observation, so that $X_{\nu+1}$
is the first post-change observation. Let $\mathsf{P}_\nu$ and
$\mathsf{E}_\nu$ denote probability and expectation when the change occurs
at $\nu+1$ for a fixed $0 \leq\nu< \infty$, and
let $\mathsf{P}_\infty$ and $\mathsf{E}_\infty$
denote the same when $\nu=\infty$ (i.e., there never is a change). A
sequential
change detection procedure is a stopping time $T$ adapted to the observations
$X_1,X_2,\dots,$ that is, $\{T \leq n\} \in{\mathscr{F}}_n$, where
${\mathscr{F}}
_n=\sigma(X_1,\dots,X_n)$ is
the sigma-algebra generated by the first $n$ observations.

Common operating characteristics of a sequential detection procedure
are the Average Run Length (ARL) to False Alarm, that is, the expected
number of observations to an alarm assuming that there is no change,
and the
Average Delay to Detection, that is, the expected delay between a
change and its
detection. The goal is to find a detection procedure that minimizes the
average detection delay subject to a bound on the ARL to false alarm.

In this paper, we will be interested in the simple changepoint problem
setting, where the
observations are independent, i.i.d. pre-change with density $f_\infty
$ and i.i.d. post-change with density $f_0$.
In other words, it is assumed that $X_n$ has density $f_\infty$ for $n
\leq\nu$ and density $f_0$
for $n > \nu$, where both $f_\infty$ and $f_0$ are known but the
changepoint $\nu$ is unknown. Therefore, the conditional density of
the sample $(X_1,\dots,X_n)$ for the fixed changepoint is
\[
p(X_1,\dots,X_n |\nu=k) = \prod_{i=1}^{k} f_\infty(X_i) \times
\prod_{i=k+1}^{n} f_0(X_k),
\]
where $\prod_{i=j}^m f(X_i) =1$ when $ j >m$.

In 1961, for detecting a change in the drift of a Brownian motion, Shiryaev
introduced a change detection procedure, which is now usually referred
to as the Shiryaev--Roberts (SR) procedure [\citeauthor{Shiryaev61}
(\citeyear{Shiryaev61}), (\citeyear{ShiryaevTPA63});
\citet{Roberts66}]. The SR procedure calls
for stopping and raising an alarm at
\begin{equation} \label{SRst}
T_{\mathrm{sr}}(A)= \inf\{n \geq1\dvtx R_n \geq A\}, \qquad\inf\{
\varnothing\}
=\infty,
\end{equation}
where
\begin{equation} \label{SRstat}
R_n = \sum_{k=0}^{n-1} \frac{p(X_1,\dots,X_n|\nu=k)}{p(X_1,\dots
,X_n|\nu=\infty)}=
\sum_{k=1}^{n} \prod_{i=k}^n \frac{f_0(X_i)}{f_\infty(X_i)}
\end{equation}
is the SR statistic, and $A>0$ is a threshold that controls the false
alarm rate.

This procedure has a number of interesting optimality properties.
In particular, if $A=A_\gamma$ is such that $\mathsf{E}_\infty
[T_{\mathrm{sr}}
(A_\gamma)]=\gamma$, then it minimizes the {\em integral average
detection delay}
\begin{eqnarray*}
\mathcal{I}(T) = \frac{\sum_{\nu=0}^\infty\mathsf{E}_\nu(T-\nu
)^+}{\mathsf{E}_\infty T}
\end{eqnarray*}
over all stopping times $T$ that satisfy
\begin{equation} \label{ARLconstraint}
\mathsf{E}_\infty T \geq\gamma,
\end{equation}
where $\gamma>1$ is a value set before the surveillance begins [cf.
\citet{PollakTartakovskySS09} and also \citet{FeinbergShiryaev07}
for the Brownian motion model].

Note that the SR statistic \eqref{SRstat} can be written recursively as
\begin{equation} \label{SRstatrec}
R_n = (1+R_{n-1}) \Lambda_n, \qquad n \geq1,\qquad  R_0 =0,
\end{equation}
where $\Lambda_n = f_0(X_n)/f_\infty(X_n)$ is the likelihood ratio.
Therefore, the classical SR statistic starts from 0.

\citet{PollakAS85} introduced a natural worst-case detection delay
measure---{\em supremum average delay to detection}
\[
\mathcal{J}_{\mathrm{P}}(T) = \sup_{0 \leq\nu<\infty} \mathbf
{E}_\nu(T-\nu|T > \nu) ,
\]
and attempted to find an optimal procedure that would minimize
$\mathcal{J}_{\mathrm{P}}
(T)$ over all procedures subject to constraint \eqref{ARLconstraint}.
Pollak's idea was to modify the SR statistic by randomization of the
initial condition $R_0$ in \eqref{SRstatrec} in order to make it an
equalizer [i.e., to make the conditional average detection delay
$\mathsf{E}
_\nu(T-\nu|T > \nu)$ independent of the changepoint $\nu$].
Pollak's version of the SR procedure starts from a random point sampled
from the quasi-stationary distribution of the SR statistic $R_n$. He
proved that this ``randomized'' procedure is asymptotically (as $\gamma
\to\infty$) optimal within an additive term of order
$o(1)$ in the sense of minimizing the supremum average detection delay
$\mathcal{J}_{\mathrm{P}}(T)$.

To be specific, let, for $B >0$,
\[
\mathsf{Q}_B(x) = \lim_{n \to\infty} \mathsf{P}_\infty\bigl(R_n \leq x
| T_{\mathrm{sr}}(B) > n\bigr)
\]
denote the quasi-stationary distribution of the SR statistic, and let
$R_n^{\mathsf{Q}_B}$ be given recursively
\begin{equation} \label{SRPstatrec}
R_n^{\mathsf{Q}_B} = (1+R_{n-1}^{\mathsf{Q}_B}) \Lambda_n,\qquad n \geq
1, \qquad
R_0^{\mathsf{Q}_B} \sim\mathsf{Q}_B,
\end{equation}
where $R_0^{\mathsf{Q}_B} \sim\mathsf{Q}_B$ means that $R_0^{\mathbf
{Q}_B}$ is a random
variable distributed according to the quasi-stationary distribution
$\mathsf{Q}_B$. The corresponding stopping time is given by
\begin{equation} \label{SRPst}
T_{\mathrm{srp}}(B)= \inf\{n \geq1\dvtx R_n^{\mathsf{Q}_B} \geq B\},
\qquad\inf\{
\varnothing\} =\infty.
\end{equation}
\citet{PollakAS85} proved that if $B=B_\gamma$ is selected so that
$\mathsf{E}_\infty[T_{\mathrm{srp}}(B_\gamma)] = \gamma$, then
\begin{equation} \label{Jpdif}
\mathcal{J}_{\mathrm{P}}(T_{\mathrm{srp}}(B_\gamma) ) - \inf_{\{T:
\mathsf{E}_\infty T \geq\gamma\}} \mathcal{J}_{\mathrm{P}}
(T) =o(1) \qquad\mbox{as $\gamma\to\infty$},
\end{equation}
where $o(1) \to0$ as $\gamma\to\infty$. We will call this
asymptotic optimality property {\em third-order asymptotic optimality}
as opposed to the second-order optimality when the corresponding
difference is bounded [i.e., $O(1)$] and the first-order optimality
when the ratio of the corresponding values tends to 1. Therefore, the
procedure given by \eqref{SRPstatrec} and \eqref{SRPst}, which we
will refer to as the Shiryaev--Roberts--Pollak (SRP) procedure, is
third-order asymptotically optimal for the low false alarm rate. Note
that this result is extremely strong since the difference between the
average detection delays in \eqref{Jpdif} is asymptotically small
while each of them is of order $O(\log\gamma)$ (i.e., both terms go
to infinity). It can be also deduced
from \citeauthor{PollakAS85} (\citeyear{PollakAS85,PollakAS87}) that the conventional SR procedure
is asymptotically minimax for a low false alarm rate within an additive
term of
order $O(1)$, that is, it is only second-order asymptotically optimal.

Since the SRP procedure is an equalizer, that is, $\mathcal
{J}_{\mathrm{P}}(T_{\mathrm{srp}})=\mathsf{E}_0
T_{\mathrm{srp}}= \mathsf{E}_\nu(T_{\mathrm{srp}}-\nu| T_{\mathrm
{srp}}>\nu)$ for all $\nu\geq0$, it is
tempting for one to conjecture that it may in fact be {\em strictly}
optimal for every $\gamma>1$. However, to date there is no proof or
disproof of this conjecture [see
\citet{YakirAS97} and \citet{MeiAS06}]. Recent work of~\citet
{MoustPolTarSS09} indicates that the SRP procedure may not be exactly
optimal and partially sheds light on this issue by
considering a generalization of the SR procedure that starts from a
specially designed deterministic point $r$. To emphasize the
dependence on the starting point, this procedure will be referred to
as the SR-$r$ procedure. Specifically, define the stopping time
\begin{equation} \label{SR-rst}
T_{\mathrm{sr}}^r(A)= \inf\{n \geq1\dvtx R_n^r \geq A\}, \qquad\inf\{
\varnothing
\} =\infty,
\end{equation}
where $R_n^r$ obeys the recursion
\begin{equation} \label{SR-rstatrec}
R_n^r = (1+R_{n-1}^r) \Lambda_n, \qquad n \geq1, \qquad R_0^r =r \geq0.
\end{equation}
Solving numerically integral equations for performance metrics for two
examples that involve Gaussian and exponential models, \citet
{MoustPolTarSS09} found that the SR-$r$ procedure (with a certain
$r=r_\gamma$ that depends on $\gamma$) uniformly outperforms the SRP
procedure, that is, $\mathsf{E}_\nu(T_{\mathrm{sr}}^r - \nu|
T_{\mathrm{sr}}^r > \nu) < \mathsf{E}_0
T_{\mathrm{srp}}$ for all $\nu\geq0$. We believe that these results present
serious evidence against optimality of the SRP procedure. However, this
may not be completely convincing since a small numerical error can be
present in such experiments.

In the present paper, we construct a counterexample where all
computations can be performed analytically. This example proves that
the SRP procedure is not optimal while the SR-$r$ procedure with a
deterministic initialization is optimal. This result answers a
long-standing question on optimality of the SRP procedure and opens a
new direction in the quest for the unknown optimum.

\section{The main theorem and integral equations for operating
characteristics} \label{s:IntEq}

Theorem~\ref{Th1} below provides a lower bound for the infimum of
Pollak's worst-case measure $\mathcal{J}_{\mathrm{P}}(T)$ which will
be used to find the
optimal changepoint detection procedure in Section~\ref{s:Ex}. Note
that a proof sketch of this lower bound has been previously given in
\citet{MoustPolTarSS09}. Here we provide a complete proof.

We first need the following lemma which establishes optimality of the
SR-$r$ procedure with respect to the integral average detection delay.

\begin{lemma}\label{Lem1}
Let
\begin{equation}\label{IRADDdef}
\mathcal{I}_r(T)= \frac{r\mathsf{E}_0T+\sum_{\nu=0}^\infty\mathbf
{E}_\nu(T-\nu
)^+}{r+\mathsf{E}_\infty T}\vadjust{\goodbreak}
\end{equation}
and let $T_{\mathrm{sr}}^r(A_\gamma)$ be the SR-$r$ detection
procedure with
$\mathsf{E}_\infty[T_{\mathrm{sr}}^r(A_\gamma)]= \gamma$. For any
$r \geq0$, the
SR-$r$ procedure minimizes $\mathcal{I}_r(T)$ over all procedures with
$\mathsf{E}
_\infty T \geq\gamma$, that is,
\begin{equation}\label{IRADDopt}
\inf_{\{T: \mathsf{E}_\infty T \geq\gamma\}}\mathcal
{I}_r(T)=\mathcal{I}_r(T_{\mathrm{sr}}
^r(A_\gamma)) .
\end{equation}
\end{lemma}

\begin{pf}
The proof of \eqref{IRADDopt} for $r=0$ is given in \citet{PollakTartakovskySS09}, Theorem 1
and Corollary 1. We now give its extension for
an arbitrary positive $r$.

Consider the following Bayesian problem, which will be denoted by
$\mathcal{B}(\pi,p,c)$. Suppose $\nu$ is a random variable
(independent of
the observations) with a zero modified geometric distribution
\begin{eqnarray*} 
P(\nu<0)=\pi, \qquad P(\nu=k)=(1-\pi)p(1-p)^{k}, \qquad k \geq0,
\end{eqnarray*}
and the losses associated with stopping at time $T$ are 1 if $T \leq
\nu
$ and $c \cdot
(T-\nu)$ if $T > \nu$, where $0\leq\pi<1$, $0<p <1$ and $c>0$ are
fixed constants. For
$\mathcal{A}\in{\mathscr{F}}$, define the probability
\[
\mathsf{P}(\mathcal{A})=\pi\mathsf{P}_0(\mathcal{A}) + (1-\pi
)\sum
_{k=0}^\infty p(1-p)^{k}\mathsf{P}_k(\mathcal{A}),
\]
and let $\mathsf{E}$ denote the corresponding expectation.

Solving $\mathcal{B}(\pi,p,c)$ requires minimization of the expected loss
\begin{eqnarray*}
\mathcal{L}_{\pi,p,c}(T)= \mathsf{P}(T\leq\nu)+ c \mathbf
{E}(T-\nu)^+
\end{eqnarray*}
or, equivalently, maximization of the expected ``gain'' $\tfrac
{1}{p}[1-\mathcal{L}_{\pi,p,c}(T)]$,
and the Bayes rule for this problem is given by the Shiryaev procedure [cf.
\citet{ShiryaevTPA63}]
\begin{eqnarray*}
T_{\pi,p,c}= \inf\bigl\{n \geq1\dvtx R_{n}^{(\pi,p)} \geq A_{\pi,p,c}\bigr\},
\end{eqnarray*}
where $A_{\pi,p,c}>\pi/(1-\pi)p$ is an appropriate threshold and
\[
R_{n}^{(\pi,p)} = \frac{\pi}{(1-\pi)p} \prod_{i=1}^n \biggl(\frac
{\Lambda_i}{1-p}\biggr)+\sum_{k=1}^n \prod_{i=k}^n \biggl(\frac{\Lambda_i}{1-p}\biggr).
\]

Let $\pi=r p$. Then, obviously, $R_{n}^{(\pi,p)} \mathop{\hbox to 19pt{\rightarrowfill}}\limits_{p\to
0} R_n^r$.

Now, it follows from \citet{PollakAS85} that there are a constant
$0<c^*<\infty$ and a
sequence $\{p_i,c_i\}_{i \geq1}$ with $p_i\to0$,
$c_i\to c^*$ as $i \to\infty$, such that $T_{\mathrm{sr}}^r(A_\gamma
)$ is the
limit of the Bayes stopping times
$T_{\pi=r p_i,p_i,c_i}$ as $i\to\infty$ and
\begin{equation} \label{LimSup}
\limsup_{p \to0, c\to c^*}
\frac{1-\mathcal{L}_{\pi=rp,p,c}(T_{\pi=rp,p,c})}{1-\mathcal
{L}_{\pi=rp,p,c}(T_{\mathrm{sr}}
^r(A_\gamma))} =1.
\end{equation}

Next, for any stopping time $T$,
\begin{eqnarray*}
\frac{\mathsf{E}(T-\nu)^+}{p} & =& \frac{\pi+(1-\pi)p}{p} \mathbf
{E}_0 T +
\frac{1-\pi}{p}\sum_{k=1}^\infty p(1-p)^{k} \mathsf{E}_k(T-k)^+
\\
& =& [r+(1-rp)] \mathsf{E}_0 T + (1-rp) \sum_{k=1}^\infty(1-p)^{k}
\mathsf{E}_k(T-k)^+
\\
& \mathop{\hbox to 19pt{\rightarrowfill}}\limits_{p \to0}& r \mathsf{E}_0 T + \sum_{k=0}^\infty
\mathsf{E}_k(T-k)^+
\end{eqnarray*}
and
\begin{eqnarray*}
\frac{\mathsf{P}(T > \nu)}{p} & =& \frac{1}{p} \Biggl(\pi+(1-\pi)p +
(1-\pi) \sum_{k=1}^\infty p(1-p)^{k}\mathsf{P}_k (T > k) \Biggr)
\\
& =& \frac{rp +(1-rp)p}{p} + \frac{1-rp}{p} \sum_{k=1}^\infty p
(1-p)^{k}\mathsf{P}_\infty(T >k )
\\
& \mathop{\hbox to 19pt{\rightarrowfill}}\limits_{p \to0}& r + \sum_{k=0}^\infty\mathsf{P}_\infty(T
> k) =r+ \mathsf{E}
_\infty T ,
\end{eqnarray*}
where we used the fact that $\mathsf{P}_k(T>k) = \mathsf{P}_\infty
(T>k)$ since by
the definition of the stopping time the event $\{T \leq k\}$ belongs to
the $\sigma$-algebra ${\mathscr{F}}_k$ and at time instant $k$ the
distribution
is still $f_\infty$.

Since
\[
\frac{1}{p}[1-\mathcal{L}_{\pi,p,c}(T)]= \frac{\mathsf{P}(T > \nu
)}{p} - c
\frac{\mathsf{E}(T-\nu)^+}{p},
\]
it follows that if $\pi=rp$, then for any stopping time $T$ with
$\mathsf{E}
_\infty T <\infty$
\begin{eqnarray*}
\frac{1}{p}[1-\mathcal{L}_{\pi=rp,p,c}(T)]  \mathop{\hbox to 19pt{\rightarrowfill}}\limits_{p \to
0}(r+\mathsf{E}_\infty T)
- c \Biggl(r \mathsf{E}_0 T + \sum_{k=0}^\infty\mathsf{E}_k(T-k)^+\Biggr),
\end{eqnarray*}
which together with \eqref{LimSup} establishes that the SR-$r$ procedure
minimizes $\mathcal{I}_r(T)$ over all stopping times that satisfy
$\mathsf{E}_\infty
T=\gamma$.
In order to prove that~\eqref{IRADDopt} holds in the class $\{T\dvtx
\mathsf{E}
_\infty T \geq\gamma\}$ it suffices to apply
the argument identical to that used in the proof of Corollary 1 in
\citet{PollakTartakovskySS09}.
\end{pf}

\begin{theorem} \label{Th1}
Let $T_{\mathrm{sr}}^r(A)$ be defined as in \eqref{SR-rst} and let
$A=A_\gamma$
be selected so that $\mathsf{E}_\infty[T_{\mathrm{sr}}^r(A_\gamma
)]=\gamma$. Then
for every $r\geq0$
\begin{equation} \label{IntADD}
\inf_{\{T: \mathsf{E}_\infty T \geq\gamma\}} \mathcal{J}_{\mathrm
{P}}(T) \geq\frac{ r \mathsf{E}
_0[T_{\mathrm{sr}}^r(A_\gamma)] + \sum_{\nu=0}^\infty\mathbf
{E}_\nu[T_{\mathrm{sr}}
^r(A_\gamma)-\nu]^+}{r+ \mathsf{E}_\infty[T_{\mathrm
{sr}}^r(A_\gamma)]} .
\end{equation}
\end{theorem}

\begin{pf}
Note first that for any stopping time $T$
\begin{eqnarray*}
\sum_{\nu=0}^\infty\mathsf{E}_\nu(T-\nu)^+ & =& \sum_{\nu
=0}^\infty\mathsf{P}
_\nu(T>\nu) \mathsf{E}_\nu(T-\nu|T>\nu)
\\
&=& \sum_{\nu=0}^\infty\mathsf{P}_\infty(T>\nu) \mathsf{E}_\nu
(T-\nu|T>\nu) ,
\end{eqnarray*}
where again we used the fact that $\mathsf{P}_\nu(T>\nu)=\mathbf
{P}_\infty(T>\nu
)$. Since
\[
\mathcal{J}_{\mathrm{P}}(T) = \sup_{k \geq0} \mathsf{E}_k(T-k|T>k)
\geq\mathsf{E}_\nu(T-\nu|T>\nu)
\qquad\mbox{for any $\nu\geq0$}
\]
and
\begin{eqnarray*}
\mathcal{J}_{\mathrm{P}}(T) & = &\frac{\mathcal{J}_{\mathrm{P}}(T)
[ r + \sum_{\nu=0}^\infty\mathsf{P}_\infty(T >
\nu)]} {r + \sum_{\nu=0}^\infty\mathsf{P}_\infty(T > \nu)}
\\
& =& \frac{r \mathcal{J}_{\mathrm{P}}(T) + \sum_{\nu=0}^\infty
\mathcal{J}_{\mathrm{P}}(T) \mathsf{P}_\infty(T >
\nu)} {r + \sum_{\nu=0}^\infty\mathsf{P}_\infty(T > \nu)},
\end{eqnarray*}
where $\sum_{\nu=0}^\infty\mathsf{P}_\infty(T > \nu)=\mathbf
{E}_\infty T$, we
obtain that for any stopping time $T$ with finite ARL to false alarm
\begin{eqnarray*}
\mathcal{J}_{\mathrm{P}}(T) & \geq&\frac{r \mathsf{E}_0 T + \sum
_{\nu=0}^\infty\mathsf{E}_\nu
(T-\nu|T>\nu) \mathsf{P}_\infty(T > \nu)} {r + \mathsf{E}_\infty T}
\\
& =&\frac{r \mathsf{E}_0 T + \sum_{\nu=0}^\infty\mathsf{E}_\nu
(T-\nu)^+} {r +
\mathsf{E}_\infty T} .
\end{eqnarray*}
Therefore,
\begin{equation} \label{Lowerbound}
\inf_{\{T: \mathsf{E}_\infty T\geq\gamma\}}\mathcal{J}_{\mathrm
{P}}(T)\geq\inf_{\{T: \mathsf{E}
_\infty T\geq\gamma\}}\mathcal{I}_r(T) ,
\end{equation}
where $\mathcal{I}_r(T)$ is defined in \eqref{IRADDdef}.

By Lemma~\ref{Lem1}, the infimum on the right-hand side in \eqref
{Lowerbound} is attained for the SR-$r$ detection procedure $T_{\mathrm{sr}}
^r(A_\gamma)$, which completes the proof.
\end{pf}

Notice that if $r$ can be chosen so that the SR-$r$ procedure becomes
an equalizer [i.e., $\mathsf{E}_0 T_{\mathrm{sr}}^r = \mathsf{E}_\nu
(T_{\mathrm{sr}}^r -\nu| T_{\mathrm{sr}}^r
> \nu)$ for $\nu\geq0$], then it is optimal since the right-hand side
in \eqref{IntADD} is equal to
$\mathsf{E}_0 T_{\mathrm{sr}}^r$ which in turn is equal to $\sup
_{\nu\geq0} \mathsf{E}_\nu
(T_{\mathrm{sr}}^r -\nu| T_{\mathrm{sr}}^r > \nu)=\mathcal
{J}_{\mathrm{P}}(T_{\mathrm{sr}}^r)$. This observation will
be used in Section~\ref{s:Ex} for proving that the SR-$r$ procedure
with a specially designed $r=r_A$ is strictly optimal for an
exponential model.

Introduce the following notation:
\begin{eqnarray*}
\delta_\nu(r) &=&\mathsf{E}_\nu(T_{\mathrm{sr}}^r-\nu)^+,\qquad \rho
_\nu(r)=\mathsf{P}_\infty
(T_{\mathrm{sr}}^r >\nu), \qquad\nu\geq0,
\\
\phi(r)&=&\mathsf{E}_\infty T_{\mathrm{sr}}^r,\qquad \psi(r) = \sum
_{\nu=0}^\infty\mathsf{E}
_\nu(T_{\mathrm{sr}}^r-\nu)^+ ,
\end{eqnarray*}
where, obviously, $\rho_0(T_{\mathrm{sr}}^r)=1$ and $\delta_0(r) =
\mathsf{E}_0 T_{\mathrm{sr}}^r$.\vadjust{\goodbreak}

In the rest of the paper we will assume for simplicity that $\Lambda
_1$ is continuous. For $i=0,\infty$, let $F_i(x) =\mathsf{P}_i
(\Lambda
_1\leq x)$ denote the distribution functions of the likelihood ratio
under the change and no-change hypotheses.

\citeauthor{MoustPolTarCS09} (\citeyear{MoustPolTarCS09,MoustPolTarSS09}) used the Markov property of
the SR-$r$ statistic \eqref{SR-rstatrec} to obtain the following
integral equations for performance metrics:
\begin{eqnarray}\label{ARLie}
\phi(r) &=& 1+\int_0^A \phi(x) \frac{\partial}{\partial x}F_\infty
\biggl(\frac{x}{1+r}\biggr) \,dx,
\\\label{ADD0ie}
\delta_0(r) &=& 1+\int_0^A \delta_0(x) \frac{\partial}{\partial
x}F_0\biggl(\frac{x}{1+r}\biggr) \,dx,
\\\label{ADDnuie}
\delta_\nu(r)&=&\int_0^{A} \delta_{\nu-1}(x) \frac{\partial
}{\partial x}F_\infty\biggl(\frac{x}{1+r}\biggr) \,dx,\qquad \nu= 1, 2,
\dots,
\\\label{rhonuie}
\rho_\nu(r)&=&\int_0^{A}\rho_{\nu-1}(x) \frac{\partial}{\partial
x}F_\infty\biggl(\frac{x}{1+r}\biggr) \,dx,\qquad \nu= 1, 2, \dots,
\\\label{psiie}
\psi(r)&=&\delta_0(r)+\int_0^{A}\psi(r) \frac{\partial}{\partial
x}F_\infty\biggl(\frac{x}{1+r}\biggr)\, dx.
\end{eqnarray}

The conditional average delay to detection of the SR-$r$ procedure is
computed as
\[
\mathsf{E}_\nu(T_{\mathrm{sr}}^r-\nu|T_{\mathrm{sr}}^r>\nu)=\frac
{\delta_\nu(r)}{\rho_\nu
(r)}, \qquad\nu\geq0,
\]
and the lower bound as
\[
\mathcal{I}_r(T_{\mathrm{sr}}^r) = \frac{r \delta_0(r) + \psi
(r)}{r+\phi(r)}.
\]

Next, we present integral equations for the operating characteristics
of the randomized SRP
procedures \eqref{SRPstatrec} and \eqref{SRPst}. Here the most crucial
problem is the computation of the quasi-stationary
distribution $\mathsf{Q}_B(x)$ of the SR statistic. By \citet{Harris63}, Theorem
III.10.1, in the continuous case the quasi-stationary
distribution exists. Its density $q_B(x)=d\mathsf{Q}_B(x)/dx$
satisfies the
following integral equation:
\begin{equation}\label{eq:QSD-eqn-int-form}
\lambda_B q_B(x)=\int_0^B q_B(r)\frac{\partial}{\partial x}F_\infty
\biggl(\frac{x}{1+r}\biggr)\,dr
\end{equation}
[see \citet{PollakAS85}], where $\lambda_B$ is the leading
eigenvalue of the linear operator associated with the kernel
\begin{eqnarray*}
K_\infty(x,r)=\frac{\partial}{\partial x}F_\infty\biggl(\frac
{x}{1+r}\biggr),\qquad x,r\in[0,B) .
\end{eqnarray*}
Thus, $q_B(x)$ is the corresponding (left) eigenfunction. It also
satisfies the constraint
\begin{equation}
\int_0^B q_B(x) \,dx=1.
\label{eq:int_to_1}
\end{equation}
Equations \eqref{eq:QSD-eqn-int-form} and \eqref{eq:int_to_1}
uniquely define $\lambda_B$ and $q_B(x)$. The equations have unique
solutions, since $\lambda_B<1$, as follows from \citet{MoustPolTarSS09}.

Once $q_B(x)$ is available we can compute the ARL to false alarm and
the average detection delay of the SRP procedure $T_{\mathrm{srp}}$
\begin{eqnarray}\label{ARLsrp}
\mathsf{E}_\infty[T_{\mathrm{srp}}(B)] & =&\int_0^B \phi(r)q_B(r)\,dr,
\\\label{ADDsrp}
\mathsf{E}_0[T_{\mathrm{srp}}(B)] & =&\int_0^B \delta_0(r)q_B(r)\,dr.
\end{eqnarray}
We recall that the SRP procedure is an equalizer: $\mathsf{E}_\nu
(T_{\mathrm{srp}}-\nu
|T_{\mathrm{srp}}>\nu)=\mathsf{E}_0T_{\mathrm{srp}}$.

The integral equations derived above are Fredholm equations of the
second kind. Usually, they do not allow for an analytical solution and
should be solved numerically. However, in the next section, we provide
an example where analytical solutions can be obtained.

\section{An example} \label{s:Ex}

Consider the exponential model with the pre-change mean 1 and the
post-change mean $\theta^{-1}$, $\theta>1$, that is, $f_\infty
(x)=e^{-x}{\mathbh{1}_{\{x\geq0\}}}$ and $f_0(x)=\theta e^{-\theta
x}{\mathbh{1}_{\{x\geq0\}}}$. We will call this model the
$\mathcal{E}(1,\theta)$-model. In the sequel we will assume that
$\theta=2$ and the thresholds in both procedures SR-$r$ and SRP do not
exceed 2.

\begin{theorem}\label{Th2} 
Assume the $\mathcal{E}(1,2)$-model. Let in the SR-$r$ procedure
$T_{\mathrm{sr}}
^{r_A}$ the initializing value be chosen as $r_A=\sqrt{1+A}-1$ and let
the threshold $A=A_\gamma$ be selected from the transcendental equation
\begin{equation} \label{treqthre}
A + (\gamma-1) \sqrt{1+A} \log(1+A) -2 (\gamma-1) \sqrt{1+A} =0.
\end{equation}
Then, for every $1<\gamma< \gamma_0$, where $\gamma_0= (1-0.5\log
3)^{-1} \approx2.2188$, the ARL to false alarm $\mathsf{E}_\infty
[T_{\mathrm{sr}}
^{r_A}(A)] =\gamma$ and the SR-$r$ procedure is minimax, that is,
\begin{equation}\label{SADDex}
\mathcal{J}_{\mathrm{P}}(T_{\mathrm{sr}}^{r_A}) = \inf_{\{T:
\mathsf{E}_\infty T \geq\gamma\}} \mathcal{J}_{\mathrm{P}}(T) .
\end{equation}

Let, in the SRP procedure, the threshold $B=B_\gamma$ be selected as
\begin{equation}\label{ThreSRP}
B = \exp\biggl\{\frac{2(\gamma-1)}{\gamma}\biggr\} -1.
\end{equation}
Then $\mathsf{E}_\infty[T_{\mathrm{srp}}(B)] =\gamma$ and $\mathcal
{J}_{\mathrm{P}}(T_{\mathrm{srp}}(B)) > \mathcal{J}_{\mathrm{P}}
(T_{\mathrm{sr}}^{r_A}(A))$ for all $1<\gamma< \gamma_0$.
Therefore, the SRP procedure is suboptimal.
\end{theorem}

\begin{pf}
Consider first the SRP procedure. As it will become apparent later,
threshold $B=B_\gamma$ in this procedure does not exceed 2 when
$\gamma< \gamma_0$. By \eqref{eq:QSD-eqn-int-form}, for $B < 2$ the
quasi-stationary density $q_B(x)= d\mathsf{Q}_B(x)/dx$ satisfies the
integral equation
\begin{eqnarray*}
\lambda_B q_B(x) = \dfrac{1}{2}\int_0^B q_B(r) \dfrac{dr}{1+r},
\end{eqnarray*}
which due to the constraint \eqref{eq:int_to_1}
yields $\lambda_B =\frac{1}{2}\log(1+B)$ and $q_B(x)=
B^{-1}\times {\mathbh{1}_{\{x\in[0,B)\}}}$. Thus, for $B<2$ the quasi-stationary
distribution $\mathsf{Q}_B(x)=x/B$ is uniform, and, moreover, it is attained
already for $n=1$ when the very first observation becomes available.

Clearly, the $\mathsf{P}_\infty$-distribution of the SRP stopping time
$T_{\mathrm{srp}}$ is geometric with the parameter $1-\lambda_B$, so
that the
ARL to false alarm is
\begin{equation}\label{ARLsrpex}
\mathsf{E}_\infty[T_{\mathrm{srp}}(B)] = \frac{1}{1-\lambda_B} =
\frac{1} {1-
(1/2) \log(1+B)} .
\end{equation}
It follows that $\mathsf{E}_\infty[T_{\mathrm{srp}}(B)]=\gamma$
when the threshold
$B=B_\gamma$ is chosen as in~\eqref{ThreSRP} and that $B<2$ whenever
$\gamma<\gamma_0$.

By \eqref{ADDsrp}, the average detection delay of the SRP procedure is
equal to
\begin{equation} \label{ADDsrpex}
\mathsf{E}_0[ T_{\mathrm{srp}}(B)] = \frac{1}{B} \int_0^B \delta
_0(r)\, dr ,
\end{equation}
so that we need to compute the ARL to detection $\delta_0(r) = \mathsf{E}_0
T_{\mathrm{sr}}^r$ of the SR-$r$ procedure which also has to be
computed for the
evaluation of the performance of the SR-$r$ procedure itself.

Assume that $A <2$. By \eqref{ADD0ie}, we have
\begin{eqnarray*}
\delta_0(r) = 1+\frac{1}{2(1+r)^2}\int_0^A\delta_0(x) x\, dx,
\end{eqnarray*}
so that
\begin{eqnarray*}
\int_0^A\delta_0(r) r \,dr & =&\int_0^Ar \,dr+\frac{1}{2}\biggl[\int_0^A\frac
{x \,dx}{(1+x)^2}\biggr]\biggl[\int_0^A\delta_0(r) r \,dr\biggr]
\\
&=&\frac{A^2}{2}+\frac{1}{2}\biggl[\log(1+A)-\frac{A}{1+A}\biggr]\biggl[\int
_0^A\delta_0(r) r \,dr\biggr],
\end{eqnarray*}
which implies that
\begin{eqnarray*}
\int_0^Ar \delta_0(r) \,dr =A^2\biggl[\frac{A}{1+A}+2\biggl(1-\frac{1}{2}\log
(1+A)\biggr)\biggr]^{-1}.
\end{eqnarray*}
Consequently,
\begin{equation} \label{delta0}
\delta_0(r) = 1+\frac{A^2}{2(1+r)^2}\biggl[\frac{A}{1+A}+2\biggl(1-\frac
{1}{2}\log(1+A)\biggr)\biggr]^{-1}.
\end{equation}

Using \eqref{ADDsrpex} and \eqref{delta0}, we find
\begin{equation}\label{ADDsrpex1}
\ \mathsf{E}_0[T_{\mathrm{srp}}(B)]  = \bar{\delta}_0(B)
 = 1+\frac{B^2}{2(1+B)}\biggl[\frac{B}{1+B}+2\biggl(1-\frac{1}{2}\log(1+B)\biggr)\biggr]^{-1}.\hspace*{-35pt}
\end{equation}

Consider now the SR-$r$ procedure. By \eqref{ARLie}, for the ARL to
false alarm $\phi(r) =\mathsf{E}_\infty[ T_{\mathrm{sr}}^r(A)]$ we have
\begin{eqnarray*}
\phi(r) = 1+\dfrac{1}{2(1+r)}\int_0^A \phi(x)\, dx,
\end{eqnarray*}
so that
\begin{eqnarray*}
\int_0^A\phi(r)\, dr =\int_0^A dr+\frac{1}{2}\biggl[\int_0^A\frac
{dr}{1+r}\biggr] \biggl[\int_0^A\phi(x) \,dx\biggr],
\end{eqnarray*}
and therefore
\begin{eqnarray*}
\int_0^A\phi(r)\, dr = A\biggl[1-\frac{1}{2}\log(1+A)\biggr]^{-1}.
\end{eqnarray*}
Consequently,
\begin{eqnarray}\label{ARLsrex}
\phi(r) = 1+\frac{A}{2(1+r)}\biggl[1-\frac{1}{2}\log(1+A)\biggr]^{-1} .
\end{eqnarray}

Recall that for $A<2$ the statistic $R_n^r$ already kicks in the
uniform quasi-stationary distribution for $n= 1$ and any $0 \leq r <A$,
so that $T_{\mathrm{sr}}^r$ is an equalizer for $\nu\geq1$ and any
$r\in[0,A)$,
that is, $\delta_\nu(r)=\bar{\delta}_0(A)$ for all $\nu\geq1$ and
$r < A$ with $\bar{\delta}_0(A)$ given by \eqref{ADDsrpex1}. This
implies that
\begin{equation}\label{supADDsrr}
\mathcal{J}_{\mathrm{P}}(T_{\mathrm{sr}}^r) = \sup_{\nu\geq0}
\mathsf{E}_\nu(T_{\mathrm{sr}}^r - \nu|T_{\mathrm{sr}}^r >
\nu) = \max\{\bar{\delta}_0(A), \delta_0(r)\}.
\end{equation}

Let $r=r_A=\sqrt{1+A}-1$, in which case $\bar{\delta}_0(A)=\delta
_0(r_A)$, that is, for this value of the head start the
SR-$r$ procedure is an equalizer for all $\nu\geq0$. Therefore, by
Theorem~\ref{Th1} the procedure $T_{\mathrm{sr}}^{r_A}$ that starts
from the
deterministic point $r_A=\sqrt{1+A}-1$ is optimal, and \eqref{SADDex}
holds if threshold $A=A_\gamma$ is selected so that $\mathbf
{E}_\infty T_{\mathrm{sr}}
^{r_A}=\gamma$. Substituting $r=\sqrt{1+A}-1$ in \eqref{ARLsrex} and
equalizing the result to $\gamma$, yields transcendental equation
\eqref{treqthre}. It is easily verified that $A_\gamma< 2$ for
$\gamma< \gamma_0$. This completes the proof of optimality of the
SR-$r$ procedure for all $1<\gamma<\gamma_0$.

In order to show that for every given $\gamma\in(1,\gamma_0)$ the
SRP procedure is inferior it suffices to show that $\mathsf{E}_\infty
[ T_{\mathrm{sr}}
^{r_A}(A)] > \mathsf{E}_\infty[T_{\mathrm{srp}}(A)]$. By \eqref
{ARLsrex}, the ARL to
false alarm of the SR-$r$ procedure is equal to
\begin{equation}\label{ARLsrex1}
\mathsf{E}_\infty[T_{\mathrm{sr}}^{r_A}(A)]=\phi(r_A) = 1+\frac
{A}{2 \sqrt
{A+1}}\biggl[1-\frac{1}{2}\log(1+A)\biggr]^{-1} .
\end{equation}
Comparing \eqref{ARLsrex1} with \eqref{ARLsrpex}, we obtain that we
have only to show that
\[
1+\frac{A}{2\sqrt{A+1}}\biggl[1-\frac{1}{2}\log(1+A)\biggr]^{-1} > \biggl[1-\frac
{1}{2}\log(1+A)\biggr]^{-1} ,
\]
that is, that $A/\sqrt{A+1} > \log(A+1)$, which holds for any $A >0$.
Thus, we conclude that the SRP procedure is suboptimal and the proof is
complete.
\end{pf}

%
%

Let, for example, $\gamma=2$. Then, by \eqref{ThreSRP} and \eqref
{ADDsrpex1}, the threshold in the SRP procedure is equal to $B
=e-1\approx1.71828$ and the average detection delay $\mathbf
{E}_0[T_{\mathrm{srp}}(B)]
=\mathcal{J}_{\mathrm{P}}(T_{\mathrm{srp}}(B)) \approx1.33275$.

For $\gamma=2$, solving the transcendental equation \eqref{treqthre}
yields $A\approx1.66485$ and the initialization point $r_A\approx
0.63244$. By \eqref{supADDsrr}, the average detection delay of the
SR-$r$ procedure
$\mathsf{E}_0[T_{\mathrm{sr}}^{r_A}(A)]=\mathcal{J}_{\mathrm
{P}}(T_{\mathrm{sr}}^{r_A}(A)) \approx1.31622$.

\begin{figure}[t]

\includegraphics{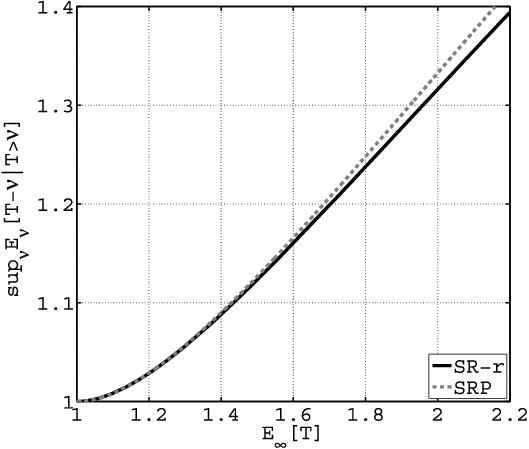}

\caption{Supremum average detection delay versus the ALR to false
alarm for $A\in(0,2)$. \label{fig}}\vspace*{6pt}
\end{figure}

Figure~\ref{fig} depicts the supremum average detection delays versus
the ARL to false alarm for the two changepoint detection procedures for
the entire range of $A \in(0,2)$.

\begin{remark*}
At an additional effort, the same conclusion can be reached in the more
general case
where the parameter of the post-change distribution $\theta>1$ and
$A,B<\theta$.
\end{remark*}

\section*{Acknowledgments}
We are grateful to George Moustakides and Moshe Pollak for useful
discussions. We would also like to thank referees and the Associate
Editor for constructive suggestions.


\printaddresses

\end{document}